\documentstyle[12pt,leqno,amscd,amssymb,verbatim]{amsart}
\newtheorem{thm}{Theorem}[section]
\newtheorem{lem}[thm]{Lemma}

\newtheorem{prop}[thm]{Proposition}

\theoremstyle{definition}

\newtheorem{defn}[thm]{Definition}

\numberwithin{equation}{thm}

\def\sB{{\cal B}}

\def\fkc{{\frak c}}
\def\fkd{{\frak d}}

\def\partial{\delta}

\def\be{\beta} 
\def\de{\delta}

\def\ph{\varphi}

\def\De{\Delta}
\def\la{\lambda}
 
\def\th{\theta} 
 
\def\La{\Lambda}

\def\tk{{\text{\tt k}}}
\def\op{{\text{\rm op}}}

\def\End{\text{\rm End}}

\def\Hom{\text{\rm Hom}}

\def\dim{\text{\rm dim\,}}

\def\im{\text{\,\rm im}}
\def\irr{{\text{\,\rm irr}}}

\begin{document}
\title[Cellular systems]%
{Finite dimensional algebras and Cellular systems}
\author{Jie Du}
\address{School of Mathematics, University of New South Wales\\
Sydney 2052, Australia}
\email{j.du@@unsw.edu.au\\
 http://www.maths.unsw.edu.au/$\sim$jied
}
\date{20/10/1999}
\thanks{Supported partially by ARC} 

\begin{abstract} We introduce the notion of a cellular system
in order to deal with quasi-hereditary algebras. We shall prove that
a necessary and sufficient condition for
an algebra to be quasi-hereditary is 
the existence of a full divisible cellular system.
As a further application, we prove that the existence
of a full local cellular system is a sufficient condition
for a standardly stratified
algebra. 
\end{abstract}

\maketitle
\begin{center}
{\it To Roger Carter on his 65th birthday}
\end{center}

\section{Introduction}

Quasi-hereditary algebras are an important class of finite dimensional
algebras with many applications to Lie theory.
Quasi-hereditary algebras have two basic formulations. The
ring-theoretic formulation defines a quasi-hereditary algebra
through the existence of a heredity chain, while in the
module-theoretic formulation it is characterized by
the notion of a highest weight category (see \cite{S},\cite{CPS},\cite{DlR}).
Motivated by the notion of cellular bases and cellular algebras
\cite{GL}, Rui and the author find a third, but incomplete,
basic formulation for a quasi-hereditary algebra $A$ over a {\it splitting}
field $\tk$ (that is, over a field $\tk$ satisfying
$\End_A(L)=\tk$ for all simple $A$-module $L$,) through the existence
of a so-called full standard basis, which is a cellular type basis
without the involution involved.
The purpose of the paper is to complete this formulation for
a quasi-hereditary algebra defined over an {\it arbitrary} field. 

A second motivation of the paper is related to a question raised
by C.C. Xi. Suppose that $A$ is an algebra over an extension field
$F$ of $\tk$ and $A$ has a cellular basis over $F$. If $F\neq \tk$,
then $A$ may have no cellular basis over $\tk$. 
(In fact, any finite dimensional algebra over a non-splitting field has no
cellular basis.) What structure can we say
for $A$ over $\tk$?

In this paper, we shall introduce the notion of a cellular system
to deal with the issues mentioned above.
Roughly speaking, a cellular system is a collection of linear injective transformations
from some $\tk$-algebras $D(\la)$ to $A$ satisfying certain 
cellular-like axioms. Various conditions imposed on $D(\la)$
will give various type of systems such as full cellular systems,
 divisible cellular systems and local cellular systems.
The main result of the paper is to prove that
an algebra is quasi-hereditary if and only if
$A$ has a full divisible cellular system.
We shall also prove that a full local cellular system
will give a standard stratification for $A$, a notion
 introduced and investigated in the recent work \cite{CPS1}
by Cline, Parshall and Scott.

We organize the paper as follows: Cellular systems and
their associated standard and costandard modules will be introduced
in \S 2 for algebras defined over a commutative ring $\tk$.
There are two useful bimodule homomorphisms $m_\la$ and $\ph_\la$
which are the main tools of our study. We will discuss them in \S3.
Full cellular systems are defined to be a system with all $\ph_\la$
surjective and will be discussed in \S 4.

In the last two sections, we shall focus on the algebras defined over 
a field $\tk$.
The main result is proved in \S5, where we first
investigate divisible cellular systems and related
representation theory.
Finally, we give some further applications to standardly stratified algebras
in \S6.

Throughout, we assume that $\tk$ is a commutative ring with 1.
By a $\tk$-algebra 
$A$ (or an algebra over $\tk$) we mean that $A$ is an
associative algebra with identity element 1, and is
finite dimensional if $\tk$ is a field.
All $A$-modules will be left modules.
A right $A$-module will be identified with a left $A^\op$-module.
So an $A^\op$-module $M^\op$ means that $M^\op$ is a right $A$-module.

\section{Algebras with cellular systems}
In this section, we assume that $\tk$ is a commutative ring with 1.

\begin{defn}\label{cel}   
 Let $A$ be a $\tk$-algebra and $\Lambda$ a poset. 
Suppose that the following conditions hold.

(a) Associated to each $\la\in\La$, there are two index sets
$I(\la),J(\la)$ and a $\tk$-algebra $D(\la)$ with 1, and,
associated to each $(i,j)\in I(\la)\times J(\la)$,
there exists a $\tk$-linear injective map 
(not necessarily an algebra homomorphism)
$c_{i,j}^\la:D(\la)\to A$ such that
$$A=\bigoplus_{\la\in\La}(\oplus_{i\in I(\la),j\in J(\la)}
c_{i,j}^\la D(\la)).$$

(b) For any $a\in A$ and $x\in D(\la)$, we have
$$\cases
(1)& a\cdot c_{i,j}^\la(x)\equiv \sum_{i'\in I(\la)}  c_{i',j}^\la
(f_{i'}^\la(a,i)x)\mod (A^{>\la})\cr 
(2)&c_{i,j}^\la(x)\cdot a\equiv \sum_{j'\in J(\la)}  
c_{i,j'}^\la(xg^\la_{j'} (j,a)) \mod(A^{>\la})\cr
\endcases$$
where $f_{i'}^\la(a,i),g^\la_{j'} (j,a)\in D(\la)$ are independent of $j$ 
and $i$, respectively,  and both independent of $x$,
and 
\begin{equation}\label{cel1}
A^{>\la}=\bigoplus_{\mu>\la}(\oplus_{i\in I(\mu),j\in J(\mu)}
c_{i,j}^\mu D(\mu)).
\end{equation}
Then,  the system $\fkc=\fkc(\La;I,J,D)=\{c_{i,j}^\la\}_{\la,i,j}$ is called a {\it cellular 
system}
of $A$ defined over the datum $(\La; I,J,D)$. 

(c) If, in addition, we assume that every $D(\la)$ is free over $\tk$, then $A$
is free. Moreover, for a given  
basis $\de^\la=\{\de_k^\la\}_{k\in K(\la)}$ 
for $D(\la)$, we
form the set 
$$\sB^\la(\de^\la)=\{c_{i,j}^\la(\de_k^\la)\mid i\in I(\la), j\in J(\la),
k\in K(\la)\}.$$
Then 
the union $\sB(\de)=\cup_{\la\in\La}\sB^\la(\de^\la)$ forms a basis for $A$.
Such a basis is called a  {\it  generalized cellular type basis} of $A$.
\end{defn}
Examples of cellular systems can be constructed from cellular or 
cellular type bases as seen from the following result.

\begin{prop}
If the $\tk$-algebra $A$ has a cellular type basis (in the sense of
\cite[(1.2.1)]{DR}) and 
$\tk_1\subseteq \tk$ is a subring (with the same identity),
then $A$ has a cellular system over $\tk_1$. 
\end{prop}

\begin{pf} Let 
$$\{a_{i,j}^\la\mid \la\in\La, (i,j)\in I(\la)\times J(\la)\}$$
be a cellular type $\tk$-basis for $A$, where $\La$ is a poset and $I(\la)$ and $J(\la)$ are index sets.
For each $\la\in\La$, let $D(\la)=\tk$ and define $\tk$-linear
injective map
$$c_{i,j}^\la:D(\la)\to A;\quad x\mapsto xa_{i,j}^\la,\forall x\in D(\la).$$
(Of course, it is $\tk_1$-linear.)
Then the condition \ref{cel}a is clearly satisfied; while
the condition \ref{cel}b follows immediately from the
corresponding condition \cite[(1.2.2)]{DR}.
Thus, we have constructed a cellular system for $\tk_1$-algebra $A$.
\end{pf}

Conversely, if $A$ has a cellular system, and suppose all $D(\la)$ are the 
same commutative ring $R$ and $A$ is also an $R$-algebra with all $c_{i,j}^\la$ $R$-linear, then the set $\{c_{i,j}^\la(1_R)\}$
is a cellular type basis for $A$ over $R$.

\begin{defn}\label{std}
Given a cellular system $\fkc(\La;I,J,D)$, there are associated
``standard'' objects $\De(\la)$ (resp. $\De^\op(\la)$) 
in the category $A$-{\bf mod} of left (resp. right) $A$-modules 
with the following properties:

(a) There are injective $\tk$-linear maps $a_i^\la:D(\la)\to\De(\la)$,
$i\in I(\la)$ (resp. $b_j^\la:D(\la)\to\De^\op(\la)$, $j\in J(\la)$)
such that 
$$\De(\la)=\bigoplus_{i\in I(\la)}a_i^\la(D(\la))
\quad (\text{ resp. } \De^\op(\la)=\bigoplus_{j\in J(\la)}b_j^\la(D(\la))).
$$

(b) The module action is given by
$$a\cdot a_i^\la(x)=\sum_{i'\in I(\la)}a_{i'}^\la(f_{i'}^\la(a,i)x)
\quad(\text{ resp. } 
b_j^\la(x)\cdot a=\sum_{j'\in I(\la)}b_{j'}^\la(xg_{j'}^\la(j,a))),$$
for all
$a\in A,  i\in I(\la), j\in J(\la), x\in D(\la)$.

Note that, since $A^{>\la}\De(\la)=0$, $\De(\la)$ is also
an $A/A^{>\la}$-module.

(c) If, in addition, we define $D(\la)$-module structure on
$\De(\la)$ (resp. $\De^\op(\la)$) naturally by
$$a_i^\la(x)\cdot y=a_i^\la(xy)\quad(\text{ resp. }
y\cdot b_j^\la(x)=b_j^\la(yx))\quad \forall x,y\in D(\la)$$
Then  $\De(\la)$ (resp. $\De^\op(\la)$) is a $A$-$D(\la)$-bimodule
(resp. $D(\la)$-$A$-bimodule) and
it is free as a $D(\la)$-module.
 
Let $\nabla(\la)=\Delta^\op(\la)^\ast$. We shall call
the modules $\Delta(\la)=\De_\fkc(\la)$ and
$\nabla(\la)=\nabla_\fkc(\la)$ the {\it standard} and
{\it costandard} objects in the category of $A$-modules {\it relative
to} the given cellular system $\fkc$.
\end{defn}

If we define $A^{\ge\la}$ similarly as in \ref{cel1} with
$<$ replaced by $\le$, we see that the quotient
$A^\la=A^{\ge \la}/A^{> \la}$ is an $A$-$A$-bimodule and
$A^\la\cong \oplus_{i,j}c_{i,j}^\la(D(\la))$ as $\tk$-module.
Note that $A^{\la}$ is an ideal of the quotient algebra
$A/A^{>\la}$.

For fixed  $(i,j)\in I(\la)\times
J(\la)$, let $\Delta(\la,j)$ (resp. $\Delta(i,\la)$) be the 
$\tk$-submodule of $A^\la$ generated by (the image of) 
$\cup_{i'\in I(\la)}c_{i',j}^\la(D(\la))$ (resp.
$\cup_{j'\in J(\la)}c_{i,j'}^\la(D(\la))$). Then,
by \ref{cel}b,  $\Delta(\la,j)$ (resp. $\Delta(i,\la)$) is 
a left (resp. right) $A$-module, and in the corresponding categories
we have $A$-module isomorphisms
$$\cases
\Delta(\la)\cong \Delta(\la,j) &\text{ via }a_{i'}^\la(x)\mapsto c_{i',j}^\la(x)+A^{>\la},\cr
\De^\op(\la)\cong \Delta(i,\la)& \text{ via }b_{j'}^\la(x)\mapsto c_{i,j'}^\la(x)+A^{>\la},\cr
\endcases
$$
for all $i'\in I(\la)$, $j'\in J(\la)$ and $x\in D(\la)$.

The following fact will be useful later on.
\begin{lem}\label{po} If $c_{i,j}^\mu(x)\De(\la)\neq0$, then $\la \ge \mu$.
\end{lem}
\begin{pf}
The hypothesis implies that some $f_{i''}^\la(c_{i,j}^\mu(x),i')$ is non-zero.
This means that $c_{i,j}^\mu(x)c_{i',j'}^\la(1)\not\equiv 0\mod A^{>\la}$.
However, $c_{i,j}^\mu(x)c_{i',j'}^\la(1)\in A^{\ge \mu}$. Therefore,
$\la\ge\mu$.
\end{pf}

\section{The homomorphisms $m_\la$ and $\ph_\la$}

We shall frequently use two bimodule homomorphisms throughout the paper.

\begin{prop}\label{m}
For any $\la\in \La$, there is an $A$-$A$-bimodule
isomorphism
$$
m_\la:\Delta(\la)\otimes_{D(\la)} \De^\op (\la)\to A^\la.
$$ 
Moreover, as left $A$-module, $A^\la$ is isomorphic to the  direct sum
$\De(\la)^{\oplus j_\la}$, where $j_\la =|J(\la)|$.
\end{prop}

\begin{pf} Consider $\tk$-bilinear map
$$f:\Delta(\la)\times \De^\op (\la)\to A^\la;\,\,
(a_i^\la(x), b_j^\la(y))\mapsto c_{i,j}^\la(xy)+A^{>\la}$$
for all $(i,j)\in I(\la)\times J(\la)$ and $x,y\in D(\la)$.
Clearly, it is $D(\la)$-balanced, i.e., $f(az,b)=f(a,zb)$
for all $a\in \De(\la)$, $b\in\De^\op(\la)$ and $z\in D(\la)$.
Therefore, $f$ induces
a $\tk$-linear homomorphism 
$$
m_\la:
\Delta(\la)\otimes_{D(\la)} \De^\op (\la)\rightarrow A^\la
$$
with $m_\la(a_i^\la(x)\otimes b_j^\la(y))=f(a_i^\la(x), b_j^\la(y))
.$
Since the restriction of $m_\la$ to $a_i^\la(D(\la))\otimes_{D(\la)}
b_j^\la(D(\la))$ is a linear isomorphism onto
$c_{i,j}^\la(D(\la))+A^{>\la}$.
So $m_\la$ is in fact a linear isomorphism.
The fact that $m_\la$ is an $A$-$A$-bimodule isomorphism follows 
from the definition \ref{std} of $\De(\la)$ and $\De^\op(\la)$.
The last assertion can be seen from the $A$-module
decomposition
$$\Delta(\la)\otimes_{D(\la)} \De^\op (\la)=
\bigoplus_{j\in J(\la)}\Delta(\la)\otimes_{D(\la)} b_j^\la(1).$$
\end{pf}

Before defining the second homomorphism, we need to look at some 
``structure constants'' with respect to the given cellular system.
 
\begin{lem}\label{mat} Let $\la\in \La$.
For any $i,i'\in I(\la)$ and $j,j'\in J(\la)$,
there exists $f^\la(j,i')\in D(\la)$ such that
$$c_{i,j}^\la(1)c_{i',j'}^\la(1)\equiv
c_{i,j'}^\la(f^\la(j,i'))\mod (A^{>\la}),$$
where $1=1_{D(\la)}$.
In general,
for any $x,y\in D(\la)$, we have 
$$c_{i,j}^\la(x)c_{i',j'}^\la(y)\equiv
c_{i,j'}^\la(xf^\la(j,i')y)\mod (A^{>\la}).$$
\end{lem}

\begin{pf}
By \ref{cel}a, we see that 
$\im (c_{i,j}^\la) \cap\im (c_{i',j'}^\la)=\{0\}$ for all $(i,j)\neq(i',j')$. So, if
 $$c_{i,j}^\la(1)c_{i',j'}^\la(1)\equiv
\sum_{i''\in I(\la)}c_{i'',j'}^\la(f^\la_{i''}(c_{i,j}^\la(1),i'))
\equiv
\sum_{j''\in J(\la)}c_{i,j''}^\la(g^\la_{j''}(j,c_{i',j'}^\la(1)))
$$ 
(mod $A^{>\la}$), then we have 
$$\sum_{i''\in I(\la)}c_{i'',j'}^\la(f^\la_{i''}(c_{i,j}^\la(1),i'))
-
\sum_{j''\in J(\la)}c_{i,j''}^\la(g^\la_{j''}(j,c_{i',j'}^\la(1)))=0,
$$
which forces
$f^\la_{i''}(c_{i,j}^\la(1),i')=g^\la_{j''}(j,c_{i',j'}^\la(1))=0$ for all $i''\neq i$, $j''\neq j'$ and
$f^\la_{i}(c_{i,j}^\la(1),i')=g^\la_{j'}(j,c_{i',j'}^\la(1))$
by the injectivity of $c_{i,j}^\la$. Putting
\begin{equation}
f^\la(j,i')=f^\la_{i}(c_{i,j}^\la(1),i')=g^\la_{j'}(j,c_{i',j'}^\la(1)),
\end{equation} 
we proved the first assertion.

Using a similar argument and noting \ref{cel}c, we have
$$c_{i,j}^\la(x)c_{i',j'}^\la(y)\equiv
c_{i,j'}^\la(f^\la_i(c_{i,j}^\la(x),i')y)
\equiv
c_{i,j'}^\la(xg^\la_{j'}(j,c_{i',j'}^\la(y)))
\mod (A^{>\la}).$$
Since $c_{i,j}^\la(x)c_{i',j'}^\la(1)\equiv
c_{i,j'}^\la(f^\la_i(c_{i,j}^\la(x),i'))\equiv
c_{i,j'}^\la(xf^\la(j,i')),$ and $c_{i,j}^\la$ is injective,
we have
 $f^\la_i(c_{i,j}^\la(x),i')=xf^\la(j,i')$. Similarly,
$g^\la_{j'}(j,c_{i',j'}^\la(y))=f^\la(j,i')y$.
Substituting gives the last relation.
\end{pf}

\begin{prop}\label{phi}  There is a $\tk$-linear map
$$\ph_{\la}: \De^\op(\la) \otimes_A \De(\la)\rightarrow D(\la)$$
such that $\ph_{\la}(b_j^\la(x)\otimes_A a_i^\la(y))=xf^{\la}(j,i)y.$
Moreover, $\ph_\la$ is a $D(\la)$-$D(\la)$-bimodule homomorphism.
Hence,   the image $\im(\ph_\la)$ is an ideal of $D(\la)$ generated
by all the $f^\la(j,i)$.
\end{prop}

\begin{pf} 
Consider the bilinear paring:
$\beta_{\la}: \Delta^\op(\la) \times \De(\la)\rightarrow D(\la)$ 
defined by
$$\beta_{\la}(b_j^\la(x), a_i^\la(y))=xf^{\la}(j,i)y.$$
We need to prove that
$$\beta_{\la}(b_j^\la(x)a, a_i^\la(y))=
\beta_{\la}(b_j^\la(x), aa_i^\la(y)),$$
for all $a\in A$.
By definition, we have
$$\text{l.h.s.}=
\beta_{\la}(\sum_{j'\in J(\la)}b_{j'}^\la(xg_{j'}^\la(j,a)), a_i^\la(y))
=x(\sum_{j'\in J(\la)}g_{j'}^\la(j,a)f^\la(j',i))y,$$
while 
$$\text{r.h.s.}=
\beta_{\la}(b_j^\la(x), \sum_{i'\in I(\la)}a_{i'}^\la(f^\la_{i'}(a,i))y))=x(\sum_{i'\in I(\la)}f^\la(j,i')f^\la_{i'}(a,i))y.$$
However, for fixed $i_0,j_0$, we have
$$\aligned
\sum_{j'\in J(\la)}g_{j'}^\la(j,a)f^\la(j',i)c_{i_0,j_0}^\la(1)
&\equiv (c_{i_0,j}^\la(1)\cdot a) c_{i,j_0}^\la(1)\cr
&\equiv c_{i_0,j}^\la(1)(a\cdot c_{i,j_0}^\la(1))\cr
&\equiv \sum_{i'\in I(\la)}f^\la(j,i')f^\la_{i'}(a,i)c_{i_0,j_0}^\la(1)\cr
&\mod (A^{>\la}).\cr
\endaligned
$$
Therefore, $\sum_{j'\in J(\la)}g_{j'}^\la(j,a)f^\la(j',i)=
\sum_{i'\in I(\la)}f^\la(j,i')f^\la_{i'}(a,i)$.
So $\be_\la$ induces a linear map $\ph_{\la}: \De^\op(\la) \otimes_A \De(\la)\rightarrow D(\la)$. It is easy to see that 
$\ph_\la$ is a $D(\la)$-$D(\la)$-bimodule homomorphism,
proving the last assertion.
\end{pf}

There is a nice relation between $m_\la$ and $\ph_\la$.
Recall from definition \ref{std} that $\De(\la)$ is an $A/A^{>\la}$-module.
 
\begin{prop}\label{ass}
 If $a,a'\in \Delta(\la)$ and $b\in \De^\op(\la) $ then 
$$m_\la(a\otimes_{D(\la)} b)a'  =a\ph_\la(b\otimes_A a'),$$ where $m_\la$
is the $A$-$A$-bimodule isomorphism given in
\ref{m}.
\end{prop}

\begin{pf}
It suffices to check the equality for
$a=a_i^\la(1),a'=a_{i'}^\la(1)$ and $b=b_j^\la(1)$.
We leave this to the reader.
\end{pf}

\section{Full cellular systems}

The homomorphism $\ph_\la$ reflects nice structure
of the algebra as we will see in the following results.

\begin{prop}\label{ful}  
Let $A$ be a $\tk$-algebra with a cellular system
$\fkc=\fkc(\La;I,J,D)=\{c_{i,j}^\la\}_{\la,i,j}$, and let $\la\in\La$.
Then the following two statements are equivalent:

(a) the ideal $A^\la$ of $A/A^{>\la}$ is idempotent;

(b) $\im(\ph_\la)=D(\la)$.
\end{prop}

\begin{pf} Let $J=A^\la$. 
Then $J^2$ is generated as $\tk$-module by all $c_{i,j}^\la(x)c_{i',j'}^\la(y)$
for all $i,j,i',j'$ and $x,y\in D(\la)$.
 Since $c_{i,j}^\la(x)c_{i',j'}^\la(y)=
c_{i,j'}^\la(xf^\la(j,i')y)$ (see \ref{mat}),  if $J^2=J$ then
 all $c_{i,j'}^\la(xf^\la(j,i')y)$, 
where $j\in J(\la)$, $i'\in I(\la)$ and $x,y\in D(\la)$, must span 
$c_{i,j'}^\la(D(\la))$. This means that 
the  ideal generated by all $f^\la(j,i')$,
which is $\im (\ph_\la)$,
must equal $D(\la)$, since $c_{i,j}^\la$ is injective. Thus, we have seen
that $J^2=J$ if and only if $\im (\ph_\la)=D(\la)$.
\end{pf}

We say that a cellular system is {\it full} if
the condition (b) in \ref{ful} holds for every $\la\in \La$. So a full system
gives idempotent ideals which in many cases have nice
homological properties (see \cite{APT}). 
Note that two classes of 
idempotent ideals --- heredity and stratifying ideals --- are
used in the study of quasi-hereditary and stratifying algebras.

If $\ph_\la$ send a `copy' of $\De^\op(\la)$ onto $D(\la)$, (Note
that this condition is stronger than the onto condition $\im(\ph_\la)=D(\la)$.)
then $\De(\la)$ will be a cyclic module.

\begin{lem}\label{Dv} For any $v\in \Delta(\la)$, if
$\ph_\la(\De^\op(\la)\otimes_A v)=D(\la)$,
then $\De(\la)=Av=A^\la v$.
\end{lem}

\begin{pf} By \ref{phi} and \ref{ass}, we have
$$
\aligned
\Delta(\la)\supseteq Av\supseteq A^\la v&=m_\la(\Delta(\la)\otimes_\tk  
\De^\op (\lambda))v\cr
&=\De(\la)\ph_\la(\De^\op (\lambda)\otimes_Av)=\Delta(\la),\cr
\endaligned$$
since $\ph_\la(\De^\op (\lambda)\otimes_Av)=D(\la)$.
Therefore,  
$\Delta(\la)=Av=A^\la v$.
\end{pf}

We can know  more properties of the modules $\De(\la)$ in this case.

\begin{prop} If $\ph_\la(\De^\op (\lambda)\otimes_Av)=D(\la)$ for some $v\in \Delta(\la)$ then

(a)  $\text{Hom}_A (\Delta(\la), \Delta(\mu))=0$
unless $\la\le  \mu$,

(b)  $\text{Hom}_A (\Delta(\la), 
\Delta(\la ))\cong D(\la)$ as $\tk$-algebra.
\end{prop}

\begin{pf} By the hypothesis and \ref{Dv}, we have
 $\Delta(\la)=Av=A^\la z$. 
Thus,  $\text{Hom}_A (\Delta(\la), 
\Delta(\mu))\neq 0$  implies that there exist  $a\in A^\la$ and $f\in \text{Hom}_A(\Delta(\la),
\Delta(\mu))$ such that
$af(v)=f(av)\neq 0$, which implies
$\la\le \mu$, by \ref{po}.   

Let $\lambda=\mu$. 
 For any $f\in \text{End}_A(\Delta(\lambda))$, $f$ is determined by
$v'=f(v)$ since $\Delta(\lambda)=A^\la v$. 
On the other hand, the hypothesis implies that
 there exists an element 
$w\in \De^\op(\la)$ such that 
$\ph_\la (w\otimes_Av)=1$.
So, by \ref{ass}, we have
$$\aligned
f(v)&=f(v\ph_\la(w\otimes_A v))=f(m_\la(v\otimes w)v)\cr
&=m_\la(v\otimes w)v'=v\ph_\la(w\otimes_Av').\cr
\endaligned
$$    
That is, every such a homomorphism is a right multiplication
by an element of $D(\la)$.
Since, for $x\in D(\la)$, right multiplication by $x$ on the
elements of $\De(\la)$ defines a left
$A$-module homomorphism on $\De(\la)$, we see that
$D(\la)=\{\ph_\lambda(w\otimes_Af(v))\mid f\in\End_A(\De(\la))\}$,
and the map  
 $\th:f\mapsto \ph_\lambda(w\otimes_Af(v))$ gives a $\tk$-linear isomorphism 
between  
$\text{End}_A(\Delta(\lambda))$ and $D(\la)$.
Now, for $f,g\in\End_A(\De(\la))$, write $f(v)=av$, for some $a\in A$.
Then, 
$$\aligned
\th(g\circ f)&=\ph_\la(w\otimes_Ag(f(v)))=\ph_\la(w\otimes_A ag(v))\cr
            &=\ph_\la(w\otimes_A av\ph_\la(w\otimes_A g(v)))\cr
         &=\ph_\la(w\otimes_A f(v))\ph_\la(w\otimes_A g(v)).\cr
\endaligned
$$
Therefore, $\th$ is an algebra isomorphism from $\End_A(\De(\la))$
to $D(\la)^\op$.
\end{pf}

\section{Divisible cellular systems}
A cellular system $\fkc(\La;I,J,D)$ of $A$ is said to be {\it divisible},
if every $D(\la)$ is a division ring. In this case,
the ring $\tk$ must be a field.
So, {\it in this section, we assume that $\tk$
is a field.}

Let  $A$ be a finite dimensional algebra over $\tk$, and
let $A$-{\bf mod} be the category of finite dimensional $A$-modules.
Recall that, for a given $A$-module
$M$, the head  hd$(M)$ of $M$ is the largest semi-simple quotient module  of $M$ and 
the radical of $M$ is the submodule $\text{rad}(M)$ of $M$ such that 
 $M/\text{rad}(M)=\text{hd}(M)$. If $L$ is a simple $A$-module, let $[M:L]$ be the multiplicity of $L$ in $M$ as a
composition factor.

We point out that the ideas behind the proofs of Theorems \ref{irr},
\ref{proj} and \ref{main} are similar to the counter-part given in 
\cite[(2.4),(4.2)]{DR}.
For completeness, we give the details of the proofs so that one may 
see some difference with the bimodule structure defined in \ref{std}c
under consideration. 

\begin{thm} \label{irr}
Let $A$ be a finite dimensional $\tk$-algebra with
a divisible cellular system $\fkc(\La;I,J,D)$ and let 
$\La_\irr=\{\la\in \La\mid\ph_\la\neq0\}$.

(a) For any $\lambda\in \Lambda_\irr$, we have 
$$\text{\rm rad} (\Delta(\la))= \{x\in \Delta(\la)\mid\ph_\la(\De^\op(\la)\otimes_A x)=0\}.$$ 
and $L(\la):=\Delta(\la)/\text{\rm rad} (\Delta(\la))$ is simple.

(b) Let $\la\in \La_\irr$.
If $L(\la)$ is a composition factor of $\Delta(\mu)$ then $\la\le \mu$ and $[\Delta(\la):L(\la)]=1$.

(c) $\{ L(\la)\mid\la\in \La_\irr\}$ is a complete set of all non-isomorphic simple
$A$-modules.
\end{thm}

\begin{pf}
Let $N_\la=\{x\in \Delta(\la)\mid\ph_\la(\De^\op(\la)\otimes_A x)=0\}$.
Then $N_\la$ is an $A$-submodule of $\De(\la)$ and  $\la\in \La_\irr$
if and only if $N_\la\neq \De(\la)$.
If $0\neq \bar v=v+N_\la\in \Delta(\lambda)/N_\la$, then
$\ph_\la(\De^\op(\la)\otimes_A v)\neq 0$, and hence
$\ph_\la(\De^\op(\la)\otimes_A v)=D(\la)$ since $D(\la)$ is a division ring
and $\ph_\la(\De^\op(\la)\otimes_A v)$ is a right ideal of $D(\la)$. 
Thus,  $\ph_\la(w\otimes_A v)=1$ for some $w$ and so
$x=x\ph_\la(w\otimes_A v)=m_\la(x\otimes w)v\in Av$ for all $x\in\De(\la)$.  
This shows that $\Delta(\la)=Av$ and $\Delta(\lambda)/N_\la=A\bar v$.
So $\Delta(\lambda)/N_\la$ is generated by its any non-zero element. 
Hence it is  a simple  $A$-module and consequently, 
$N_\la\supseteq \text{rad}(\Delta(\lambda))$.  

If  $\text {rad}(\Delta(\lambda))\neq N_\la$ then there is a simple
 module $L$ in hd$(\Delta(\la))$ such that    the projection $\eta$ from $\Delta(\la)$ onto $L$
does not map $N_\la$ to  zero. Therefore  $\eta(N_\la)=L=\eta(\De(\la))$. Thus, we have
$\eta(v)=\eta(u)$ for some $u\in N_\la$.
Since $\ph_\la (w\otimes_A u)=0$, we have  
$$
\aligned
\eta(v)=\eta(v\be_\la(w\otimes_A v))&=m_\la(v\otimes w)\eta(v)\cr
&=
\eta(m(v\otimes v)u)=\eta(v\ph_\la(w\otimes_Au))=0.\cr
\endaligned
$$
So $\eta(\Delta(\lambda))=0$, a contradiction. Therefore, rad$\,\De(\la)=N_\la$, proving (a).

We now prove (b). If $L(\la)$ is a composition factor of $\Delta(\mu)$
 then there is an 
 $A$-homomorphism $\eta: \Delta(\la)\rightarrow \Delta(\mu)/N$ for some $A$-submodule
$N\subset \Delta(\mu)$ such that $\text{im} (\eta)\cong L(\la)$. Since 
$\la\in\La_\irr$, we have $\Delta(\la)=A^\la v$ for some $v\in \Delta(\la)$.
Thus, $\eta\neq 0$, i.e., $a\eta(v)\neq0$ for some $a\in A^\la$,
 implies that there exist $c_{i,j}^\la(1)$ and $a_k^\mu(1)$ such that
$c_{i,j}^\la(1)a_k^\mu(1)\not \in N$, forcing $\la \le \mu$ by \ref{po}.
 
If $\mu=\la$, we claim that 
$\eta$ is surjective. Indeed,
suppose $\eta(v)=v'+N$. Then
$$\aligned
\eta(v)&=\eta(v\ph_\la(w\otimes_A v))=
\eta(m_\la(v\otimes w)v)\cr
&=m_\la(v\otimes w)v'+N=
v\ph_\la(w\otimes_A v')+N.\endaligned
$$
Since $\ph_\la(w\otimes_A v')\neq 0$, $\eta$ is surjective.
By the claim, we have 
$\Delta(\la)/N\cong L(\la)$, and       hence, $[\Delta(\la):L(\la)]=1$.

 To prove (c), we first note from the argument above that
$L(\la)\cong L(\mu)$ 
implies $\la=\mu$. 
We now prove that, for any simple $A$-module $L$, 
$L\cong L(\la)$ for some $\la\in\La_\irr$.
We order $\La$ linearly: $\mu=\la_1, \la_2,\ldots, \la_n$
such that 
$\la_i\ge \la_j$ implies $i\le j$, and define
$J_k=\bigoplus_{r=1}^k(\oplus_{i,j}c_{i,j}^{\la_r}(D(\la_r))$.
Then, we have a chain of ideals of $A$ 
\begin{equation}\label{seq}
0=J_0\subset J_1\subset     
   \cdots\subset J_m=A
\end{equation}
Let $f: A\rightarrow L$ be an epimorphism, 
and let $J_i$ be the minimal
 ideal in the above filtration such that
 $f|_{J_i}\neq 0$.
 Then $f$ induces an epimorphism 
$f:A^{\la_i}\cong J_i/J_{i-1}\rightarrow L$.
Since $J_{i-1}L=0$ and $J_iL\neq0$, it follows that $(A^{\la_i})^2\neq0$,
which implies $\ph_{\la_i}\neq0$. Therefore, $\la_i\in\La_\irr$.
Now, since $A^{\la_i}\cong\De(\la_i)^{\oplus j_{\la_i}}$ (\ref{m}),
we obtain an epimorphism 
$f:\Delta(\la_i)\rightarrow
L$. Therefore, by (b), $L(\la_i)\cong L$.
\end{pf}

Let
$\De=\{M\in\text{Ob\,}(A\text{-\bf mod})\mid M\cong\De(\la)\text{ for some }
\la\in\La\}.$
A $\De$-filtration of a module $M$ is a sequence of submodules:
$$0=M_0\subseteq M_1\subseteq\cdots \subseteq M_n=M$$
such that $M_i/M_{i-1}\in\De$.

\begin{thm}\label{proj}
Let $A$ be a finite dimensional $\tk$-algebra with
a divisible cellular system $\fkc(\La;I,J,D)$.
 Then the projective cover $P(\la)$ of $L(\la)$ with $\la\in\La_\irr$
has a $\De$-filtration. If $[P(\la):\De(\mu)]$ denotes the number of sections 
 isomorphic to $\De(\mu)$ in such
a filtration, then $[P(\la):\De(\mu)]\neq 0 \Rightarrow
\mu\ge\la$ and $[P(\la):\De(\la)]=1$.
\end{thm}

\begin{pf}
  Let $P=P(\la)$, and let $i$ be the minimal index in the sequence
\ref{seq} such that $J_iP= P$. Then
$J_i=M\oplus A^{\ge\la_i}$ where
 $$M=\bigoplus_{j\le i,\la_j\not>\la_i}(\oplus_{r,s}c_{r,s}^{\la_j}(D(\la_j)).$$
 is an ideal contained in $J_i$, and
 $P=J_iP=MP\oplus A^{\ge\la_i}P=A^{\ge\la_i}P$ since
$P$ is indecomposable and $MP=0$ by \ref{po}.
For any $\nu\in\La$, put $A^\nu P=A^{\ge\nu}P/A^{>\nu}P$.
Then, the above argument shows that $\la_i$ is the minimal $\nu$
(w.r.t. the partial ordering on $\La$)
 with $A^\nu P\neq0$.
In particular, $\la_i\le \la$ since we certainly have $A^\la P\neq0$.
On the other hand, since  $A^{\la_i} P\cong A^{\la_i}\otimes_A P $ is a homomorphic image of $P$. 
We have
$$ \aligned 
\text{ Hom}_A (P , \Delta(\lambda_i)) &         
\supseteq\text{Hom}_A( A^{\la_i} \otimes_A P, \Delta(\lambda_i)) 
\cr 
&\cong \text{Hom}_{D(\la_i)} (\De^\op(\la_i)\otimes_A P, 
\text{End}_A(\Delta(\la_i)) \quad \text{ by \ref{m}}\cr 
&
\cong \text{Hom}_{D(\la_i)}(\De^\op(\lambda_i)\otimes_A P, D(\la_i)) \quad 
\text{by \ref{phi}}\cr &
\neq 0,\quad\text{ since }\De^\op(\lambda_i)\otimes_A P\neq0.\endaligned $$
So $L(\la)$ is a composition factor 
of $\Delta(\la_i)$. By \ref{irr} we have $\la\le \la_i$.
Therefore, $\la_i= \la$ and $P=A^{\ge\la}P$.
Now, $A^{\ge\la}$ has a filtration with sections $A^\mu$, $\mu\ge\la$.
It follows that $P$ has a filtration with
sections isomorphic to $A^\mu\otimes_A P\cong \De(\mu)\otimes_{D(\mu)}(\De^\op(\mu)\otimes_AP)$. Therefore,
$P$ has a $\De$-filtration whose sections $\De(\mu)$ satisfies
$\mu\ge \la$. Finally,
since  $\dim_{D(\la)}(\De^\op(\la)\otimes_AP)=\dim_{D(\la)}
\Hom_{D(\la)}(\De^\op(\la)\otimes_AP,D(\la))\le\dim_{D(\la)}
\Hom_A(P(\la),\De(\la))=1$, we have
$A^{\la}\otimes_A P(\la)\cong \De(\la)$,
proving $[P(\la):\De(\la)]=1$. 
\end{pf}

We are now ready to prove the main result of the paper.

\begin{thm}\label{main} Let $A$ be a finite dimensional algebra over a field
$\tk$. Then $A$ is quasi-hereditary if and only if $A$ has 
a divisible cellular system $\fkc(\La;I,J,D)$ with $\La=\La_\irr$.
The latter is equivalent to that $A$ has 
a full divisible cellular system
\end{thm}

\begin{pf} The ``if'' part follows from \cite[Theorem 3.6]{CPS}
since, by Theorems \ref{irr} and \ref{proj}, the category $A$-{\bf mod}
is a highest weight category. It is also easy to check directly that
the sequence \ref{seq} is a heredity chain in this case.

Conversely, suppose now that $A$ is quasi-hereditary. Then,
there is a poset $\La$ on which the highest weight category
$A$-{\bf mod} is defined.
We now construct a divisible cellular system by induction on $\La$.

If $\La=\{\la\}$ has a single element, then $A$ itself is a heredity ideal.
So there exists a primitive idempotent
$e$ of $A$ such that $A\cong Ae\otimes _{eAe}eA$.
Since $Ae$ is projective indecomposable,
$eAe$ is a division ring. Thus, $Ae$
(resp. $eA$) is a right (resp. left)
vector space over $eAe$.
Choose $eAe$-bases $\{a_i\}_{i\in I}$ and $\{b_j\}_{j\in J}$  for $Ae$ and
$eA$, respectively, and define
$$c_{i,j}:eAe \to A$$
sending $x\in eAe$ to $c_{i,j}(x)=a_i^\la x\otimes_{eAe}b_j^\la$.
Putting $D=eAe$, one obtains a cellular system
$\fkc(\La;I,J,D)=\{c_{i,j}\}_{i,j}$. Clearly, it is divisible and
$\La_\irr=\La$.

Assume  $|\La|>1$ and let $\mu\in \La$ be a maximal element.   
We order $\La$ linearly: $\mu=\la_1, \la_2,\ldots, \la_n$
such that 
$\la_i\ge \la_j$ implies $i\le j$ and 
$\{\la_2,\cdots,\la_t\}=\{\nu\in\La:\nu\not\le\mu\}$.
Then,  there is a heredity chain of $A$ 
$$
0=J_0\subset J=J_1\subset \cdots\subset J_m=A,
$$
such that $J_i/J_{i-1}\cong A_if\otimes_{fA_if}fA_i$ for some (primitive)
idempotent $f$ of $A_i:=A/A_{i-1}$, and $_A(A_if)\cong\De(\la_i)$,
the standard object corresponding to $\la_i$.

Let $B=A/J$ and $\pi:A\to B$ the natural epimorphism. 
Then $B$ is a quasi-hereditary algebra 
with the weight poset $\La_1=\La\backslash\{\mu\}$. By induction, 
$B$ has a divisible cellular system
$\fkd=\fkd(\La_1,I,J,D)=\{d_{i,j}^\la\}_{\la,i,j}$
with the property
$$J_r/J=(\oplus_{i,j}
d_{i,j}^{\la_2}D(\la_2))\oplus\cdots\oplus
(\oplus_{i,j}
d_{i,j}^{\la_r}D(\la_r))
,\,\,\, \forall r\in[2,m].$$

We now inductively construct a cellular system $\{c_{i,j}^{\la_k}\}$ for $A$.
Since $J=J_1$ is a heredity ideal of $A$, there is an idempotent $e$
of $A$ such that  $J\cong Ae\otimes_{eAe} eA$ 
and $Ae\cong\De(\mu)$ is a projective indecomposable module.
So $eAe$ is a division algebra. We extend the map $D$ from
$\La_1$ to $\La$ by setting $D(\mu)=eAe$.
As in the discussion for $|\La|=1$, we may define
injective linear maps $c_{i,j}^{\mu}:D(\mu)\to J$,
where $i\in I_\mu$ and $j\in J_\mu$ for some index sets $I_\mu, J_\mu$,
 such that $J=\oplus_{i,j}c_{i,j}^\mu(D(\mu))$.
Now, the maps $I$ and $J$ are also extended by setting
$I(\mu)=I_\mu$ and $J(\mu)=J_\mu$.

Since $\la_k\not\le\mu$, $2\le k\le t$, we have
$J_t\cong J_t/J\oplus J$. So $J_t/J$ can be viewed as an ideal of $A$, and we 
may define $c_{i,j}^{\la_k}=d_{i,j}^{\la_k}$ for $2\le k\le t$.

Assume now $k>t$.
Using the linear injections
$d_{i,j}^{\la_k}:D(\la_k)\to J_k/J$, we may easily define 
linear injective maps
$c_{i,j}^{\la_k}:D(\la_k)\to J_k$ such that
$\pi\circ c_{i,j}^{\la_k}=d_{i,j}^{\la_k}$ and
the following diagram commute:
$$\matrix
D(\la_k)&\overset{c_{i,j}^{\la_k}}\longrightarrow&J_k\cr
\downarrow^{d_{i,j}^{\la_k}}&&\downarrow^{\pi_k}\cr
J_k/J&\overset{\bar\pi_k}\longrightarrow&J_k/J_{k-1}\cr
\endmatrix
$$
where $\bar\pi_k$ and $\pi$ are natural maps.
Note that, since $\bar\pi_k(d_{i,j}^{\la_k}(D(\la_k))=J_k/J_{k-1}$, we have
$\pi_k(c_{i,j}^{\la_k}(D(\la_k))=J_k/J_{k-1}$.
Therefore, one sees easily that
$A=\oplus_{\la,i,j} c_{i,j}^\la(D(\la)) $ as vector space.

It remains to prove that the system $\{c_{i,j}^\la\}_{\la,i,j}$ satisfies
the conditions in \ref{cel}b. The conditions hold clearly for $\la\in\{
\la_1,\cdots,\la_t\}$.
Since $\{d_{i,j}^\la\}$ is a cellular system of $B$, we have
for $\la=\la_k$ with $k>t$, $a\in A$ and $x\in D(\la)$, 
$$\aligned
\pi(a)\cdot d_{i,j}^\la(x)&\equiv\sum_{i'\in I(\la)}  d_{i',j}^\la
(f_{i'}^\la (\pi(a),i)x) 
\mod (B^{>\la})\cr
&\equiv\sum_{i'\in I(\la)}  d_{i',j}^\la
(f_{i'}^\la (\pi(a),i)x) 
\mod(J_{k-1}/J).\cr
\endaligned$$ 
Thus, from the definition of $c_{i,j}^\la$ and noting
$J\subset A^{>\la}$, we have
$$\aligned
a\cdot c_{i,j}^\la(x)&\equiv\sum_{i'\in I(\la)}  c_{i',j}^\la
(f_{i'}^\la (\pi(a),i)x) 
\mod(J_{k-1})\cr
&\equiv\sum_{i'\in I(\la)}  c_{i',j}^\la
(f_{i'}^\la (\pi(a),i)x) 
\mod (A^{>\la}),\cr
\endaligned$$
where the $f_{i'}^\la (\pi(a),i)$ is independent of $j$ and $x$.
By a symmetric argument, we see that a
 similar relation holds for $c_{i,j}^\la(x)\cdot a$.
Therefore, we have obtained a cellular system which is clearly
divisible
and $\La=\La_\irr$.
\end{pf}

\section{Local cellular systems and standard stratifications}

We shall assume again in this section that $\tk$ is a field
and $A$ is finite dimensional over $A$.

A cellular system $\fkc(\La;I,J,D)$ of $A$ is said to be {\it local},
if every $D(\la)$ is a local ring.

Recall from \cite[(2.1)]{CPS1} that an ideal of $A$ is called a 
(left) {\it standard
stratifying ideal} if $J=AeA$ for some idempotent $e\in A$ and $_AJ$ is a 
projective $A$-module. The algebra $A$ is said to be (left) {\it standardly
stratified} of length $n$ if $A$ has a chain
$$0=J_0\subset J_1\subset     
   \cdots\subset J_n=A$$
of ideals such that $J_i/J_{i-1}$ is a standard stratifying ideal in
$A/J_{i-1}$ ($i=1,2,\cdots,n$). In this case, the chain
is called a (left) {\it standard stratification} of $A$.
If the condition that $_AJ$ is projective is replaced
by that $J_A$ is projective, then we obtain the notion of
a right standard stratifying ideal, etc..

\begin{thm} Let $A$ be a finite dimensional $\tk$-algebra. 
Suppose that $A$ has a local cellular system $\fkc(\La;I,I,D)$ which is full,
i.e., $\im(\ph_\la)=D(\la)$ for all $\la\in\La$.
Then
$A$ is (left and right) standardly stratified of length $|\La|$.
\end{thm}

\begin{pf} We order $\La$ linearly: $\la_1, \la_2,\ldots, \la_n$
such that 
$\la_i\ge \la_j$ implies $i\le j$ and define
$$J_k=\bigoplus_{r=1}^k(\oplus_{i\in I(\la),j\in J(\la)}c_{i,j}^{\la_r}D(\la_r)).$$
Then, $A^{\ge\la_k}\subseteq J_k$, and we obtain a chain of ideals
\begin{equation}\label{str}
0=J_0\subset J_1\subset \cdots\subset J_m=A.
\end{equation}
We now prove that this chain is a standard stratification.

Since each $D(\la)$ is local and $\im(\ph_\la)=D(\la)$, 
 there exists at least one
$z=f^\la(j,i)$, which is not in the unique
maximal ideal of $D(\la)$. Hence, $z$ is invertible.
Thus,
$$c_{i,j}^\la(z^{-1})^2\equiv c_{i,j}^\la(z^{-1}f^\la(j,i)z^{-1})
\equiv c_{i,j}^\la(z^{-1})\mod(A^{>\la}).$$
Let $e$ be the image of $c_{i,j}^\la(z^{-1})$ in $\bar A=A/J_{k-1}$, where
$\la=\la_{k}$. Then
$e$ is an idempotent of $A$ and
$ J_k/J_{k-1}= \bar Ae\bar A$ 
has a basis consisting of the images of the elements of the form
$c_{i',j'}^\la(x)=c_{i',j}^\la(x)ec_{i,j'}^\la(z^{-1})$.
Therefore, $\bar Ae\cong\De(\la)$ and
$J_k/J_{k-1}$ is  a (left and right) projective 
$\bar A$-module by \ref{m}. Consequently,
it is a standard stratifying ideal and the chain \ref{str}
is  a standard stratification of $A$.
\end{pf}

We remark that, by the theorem,
each $\De(\la)$ is a projective indecomposable  
$A/A^{>\la}$-module, but its simple head may have multiplicity
more than 1. 

From the above proof, we have immediately the 
following result.

\begin{thm} Let $A$ be a finite dimensional $\tk$-algebra. 
Suppose that $A$ has a cellular system $\fkc(\La;I,I,D)$ with the 
following property: For any $\la\in\La$, there exists
$(i,j)\in I(\la)\times J(\la)$ such that $f^\la(j,i)$
is invertible in $D(\la)$.
Then
$A$ is (left and right) standardly stratified.
\end{thm}

\bigskip

\noindent{\small
{\bf Acknowledgment.} The author would like to thank C.C. Xi
for several helpful discussions. 
 The paper was written while the author 
was on leave from the University of New South Wales at the Universities
of Aarhus and Bielefeld. He wishes to thank these two universities
for their hospitality during the writing of the paper.}

\end{document}